\documentstyle[12pt]{article}

\newcommand{\be}{\begin{equation}}
\newcommand{\ee}{\end{equation}}
\newcommand{\ba}{\begin{eqnarray}}
\newcommand{\ea}{\end{eqnarray}}
\newcommand{\ban}{\begin{eqnarray*}}
\newcommand{\ean}{\end{eqnarray*}}
\newcommand{\n}{\nonumber \\}

\newcommand{\eq}[1]{(\ref{#1})}
\newcommand{\sfrac}[2]{{\textstyle \frac{#1}{#2}}}
\newcommand{\ignore}[1]{}
\newcommand{\no}{{\textstyle{\circ\atop\circ}}}

\font\csc=cmcsc10 scaled\magstep1
\font\tennn=msbm10
\font\twelvenn=msbm10 scaled\magstep1

\newcommand{\Sbm}[1]{\leavevmode\raise-.15ex\hbox{\twelvenn #1}}
\newcommand{\sbm}[1]{\leavevmode\raise-.15ex\hbox{\tennn #1}}
\newcommand{\bC}{\Sbm{C}}

\newcommand{\bZ}{\Sbm{Z}}
\newcommand{\bz}{\sbm{Z}}

\begin{document}

\renewcommand{\thefootnote}{\fnsymbol{footnote}}
\font\csc=cmcsc10 scaled\magstep1

{\baselineskip=14pt
 \rightline{
 \vbox{
       \hbox{DPSU-96-16}
       \hbox{UT-762}
       \hbox{November 1996}
}}}

\vskip 11mm
\begin{center}

{\large\bf Quantum Deformation of the $W_N$ Algebra\footnote[2]{
Talk presented by S.O. 
at the Nankai-CRM joint meeting on ``Extended and Quantum Algebras
and their Applications to Physics", Tianjin, China, August 19-24, 1996. 
To appear in the CRM series in mathematical physics, Springer Verlag.
}}

\vspace{15mm}

{\csc Hidetoshi AWATA}\footnote[1]{JSPS fellow}$^{1}$,
{\csc Harunobu KUBO}$^{*2}$,
\\ 
\vskip.05in
{\csc Satoru ODAKE}$\,{}^3$ and
{\csc Jun'ichi SHIRAISHI}$\,{}^4$
{\baselineskip=15pt
\it\vskip.35in 
\setcounter{footnote}{0}\renewcommand{\thefootnote}{\arabic{footnote}}
\footnote{e-mail address : awata@rainbow.uchicago.edu}
James Frank Institute and Enrico Fermi Institute,
University of Chicago,\\
5640 S. Ellis Ave., Chicago, IL 60637, U.S.A.
\vskip.1in 
\footnote{e-mail address : kubo@hep-th.phys.s.u-tokyo.ac.jp}
Department of Physics, Faculty of Science \\
University of Tokyo, Tokyo 113, Japan \\
\vskip.1in 
\footnote{e-mail address : odake@azusa.shinshu-u.ac.jp
}
Department of Physics, Faculty of Science \\
Shinshu University, Matsumoto 390, Japan\\
\vskip.1in 
\footnote{e-mail address : shiraish@momo.issp.u-tokyo.ac.jp}
Institute for Solid State Physics, \\
University of Tokyo, Tokyo 106, Japan \\
}
\end{center}

\vspace{7mm}
\begin{abstract}
We review the $W_N$ algebra and its quantum deformation,
based on free field realizations.
The (quantum deformed) $W_N$ algebra is defined through
the (quantum deformed) Miura transformation,
and its singular vectors realize the Jack (Macdonald) polynomials.
\end{abstract}

\vspace{2mm}
q-alg/9612001

\vfill\eject
\setcounter{footnote}{0}\renewcommand{\thefootnote}{\arabic{footnote}}

\section{Introduction}

The conformal field theory (CFT)\cite{rBPZ}, whose examples are string
theories as a world sheet theory and statistical critical phenomena,
has made remarkable progress contacting with various branches of mathematics.
CFT is a massless theory and its infinitely many conserved quantities are
controlled by the its symmetry algebra, the Virasoro algebra.
By perturbing CFT, it becomes a massive theory and there does not exist
the Virasoro algebra any longer. However, if we add a good perturbation,
the theory has still infinitely many conserved quantities and they are called
massive integrable theory (MIT).
Behind the conserved quantities, there exist symmetries.
So we would like to clarify
$$
  \mbox{What symmetry ensures the integrability of MIT?}
$$
In some cases the quantum group, Yangian, degenerate affine Hecke algebra, 
etc. have an important role. For example\cite{rDFJMN}, in spin one-half XXZ
spin chain, its correlation functions were derived by using the quantum
affine Lie algebra symmetry $U_q(\widehat{sl}_2)$.
However these symmetries correspond to the current algebra (affine Lie
algebra) symmetry in CFT, not to the Virasoro algebra.
Naively, ``quantum deformation" ($q$-deformation) of the Virasoro algebra
has been expected. After the name of the quantum group had been known to
physicists, several attempts to construct $q$-deformed Virasoro algebra
have been made. But satisfactory $q$-deformation of the Virasoro algebra 
had not been obtained(at least to the talker)\footnote{
Generally speaking, what is $q$-deformation?
Although there is no precise definition of $q$-deformation, we would like
to define $q$-deformation in the following way;
(\romannumeral1) Theory deformed by adding one parameter $q$, which reduces
to the original theory in the $q\rightarrow 1$ limit,
(\romannumeral2) (Some) Properties of the original theory remains in
the $q$-world.
The condition (\romannumeral2) is obscure and arbitrary.
So $q$-deformation is not unique. There exist ``good" $q$-deformation,
``bad" $q$-deformation and ``usual" $q$-deformation. 
}, because of the lack of definite guiding principle.

Last summer reasonable $q$-deformation of the Virasoro and $W_N$ algebras
was obtained\cite{rSKAO,rAKOS,rFF} (see also important works \cite{rFR,rLP1}).
Here we will review these $q$-deformation.
Putting aside the above physical motivation, we take the following point
of view,
$$
  \mbox{Algebra, Representation Theory, Free Field Realization}.
$$
Our guiding principle is the following.
First we note the two facts:
\begin{enumerate}
\item
In the free field realization, the singular vectors of the Virasoro and $W_N$ 
algebras realize the Jack symmetric polynomials\cite{rMY,rAMOS}.
\item
The Jack symmetric polynomials have the good $q$-deformation, the Macdonald
symmetric polynomials\cite{rM}.
\end{enumerate}
Based on these, setting up the following question seems to be natural;
\begin{itemize}
\item
Construct the algebras whose singular vectors in the free field realization
realize the Macdonald symmetric polynomials.
\end{itemize}
The resultant algebra are worth being called quantum deformation 
($q$-deformation) of the Virasoro and $W_N$ algebras in this sense.
We call these as $q$-Virasoro and $q$-$W_N$ algebras.
We illustrate this scenario by a figure,

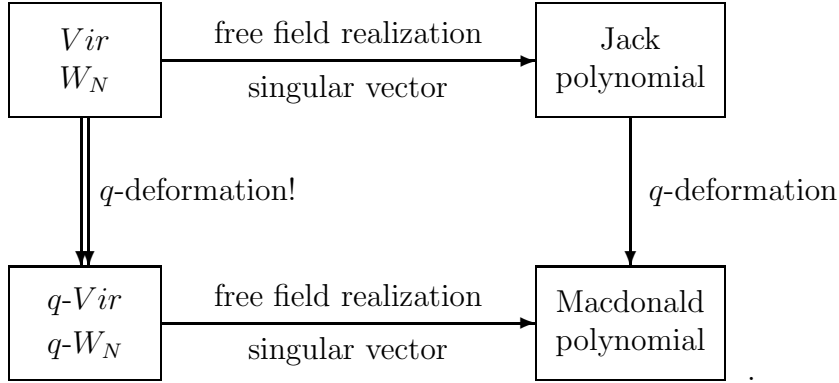
\setlength{\unitlength}{1mm}
\begin{picture}(140,65)(-10,3)
\put(10,45){\framebox(20,15){\shortstack{$Vir$ \\ \\ $W_N$}}}
\put(10,10){\framebox(20,15){\shortstack{$q$-$Vir$ \\ \\ $q$-$W_N$}}}
\put(80,45){\framebox(25,15){\shortstack{Jack \\ \\ polynomial}}}
\put(80,10){\framebox(25,15){\shortstack{Macdonald \\ \\ polynomial}}}
\thicklines
\put(30,52.5){\vector(1,0){50}}
\put(30,17.5){\vector(1,0){50}}
\put(19.5,45){\vector(0,-1){20}}
\put(20.5,45){\vector(0,-1){20}}
\put(92.5,45){\vector(0,-1){20}}
\put(30,52.5){\makebox(50,7.5){free field realization}}
\put(30,45){\makebox(50,7.5){singular vector}}
\put(30,17.5){\makebox(50,7.5){free field realization}}
\put(30,10){\makebox(50,7.5){singular vector}}
\put(25,22.5){\makebox(20,25){$q$-deformation!}}
\put(97.5,22.5){\makebox(20,25){$q$-deformation}}
\put(108,10){.}
\end{picture}

\noindent
A priori these algebras have nothing to do with the symmetry of massive
integrable models. But we believe that they are related very closely.
In fact their relations have been revealed
gradually\cite{rL,rFr,rLP2,rKa,rAKMOS,rAJMP,rJLMP,rJKM}.

In this talk we would like to review the $W_N$ and $q$-$W_N$ algebras
based on free field realizations. In section 2 the $W_N$ algebra is defined
by the Miura transformation, and singular vectors realize the Jack
symmetric polynomials. In section 3 $q$-$W_N$ algebra is introduced in a 
similar manner. We show that singular vectors realize the
Macdonald symmetric polynomials and $q$-$W_N$ reduces to $W_N$ 
in the $q\rightarrow 1$ limit. Section 4 is devoted to the discussion.
In appendix we present explicit examples.

Before going section 2, we recall the parameters of the Jack and Macdonald
symmetric polynomials (for precise definitions of these polynomials, see
other talks in this meeting).
The Jack symmetric polynomials $J_{\lambda}(x;\beta)$
have one parameter $\beta$, and
the Macdonald symmetric polynomials $P_{\lambda}(x;q,t)$
have two parameters $q$ and $t$.
The relation among the Macdonald, Jack, Hall-Littlewood and
Schur symmetric polynomials is\cite{rM},

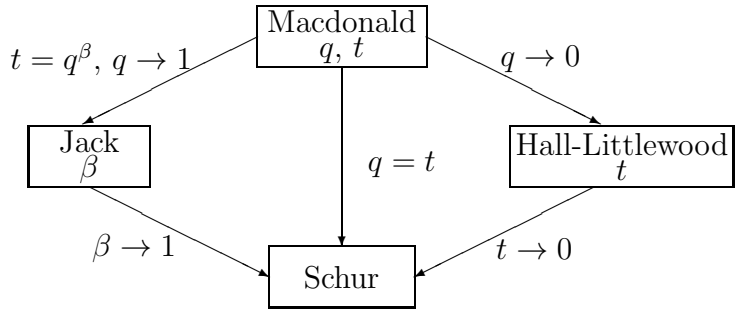
\setlength{\unitlength}{0.8mm}
\begin{picture}(140,58)(-20,8)
\put(48,50){\framebox(28,10){$\mbox{Macdonald}\atop\mbox{$q$, $t$}$}}
\put(50,10){\framebox(24,10){Schur}}
\put(10,30){\framebox(20,10){$\mbox{Jack}\atop\mbox{$\beta$}$}}
\put(90,30){\framebox(37,10){$\mbox{Hall-Littlewood}\atop\mbox{$t$}$}}
%
\put(48,55){\vector(-2,-1){29}}
\put(76,55){\vector(2,-1){29}}
\put(20,30){\vector(2,-1){30}}
\put(104,30){\vector(-2,-1){30}}
\put(62,50){\vector(0,-1){30}}
\put(12,47.5){\makebox(20,7.5){$t=q^{\beta}$, $q\rightarrow 1$}}
\put(85,47.5){\makebox(20,7.5){$q\rightarrow 0$}}
\put(17.5,16){\makebox(20,7.5){$\beta\rightarrow 1$}}
\put(84,16){\makebox(20,7.5){$t\rightarrow 0$}}
\put(62,30){\makebox(20,7.5){$q=t$}}
\put(130,30){.}
\end{picture}
%
\ignore{
\setlength{\unitlength}{1mm}
\begin{picture}(140,65)(0,4)
\put(48,50){\framebox(24,10){$\mbox{Macdonald}\atop\mbox{$q$, $t$}$}}
\put(50,10){\framebox(20,10){Schur}}
\put(10,30){\framebox(20,10){$\mbox{Jack}\atop\mbox{$\beta$}$}}
\put(90,30){\framebox(30,10){$\mbox{Hall-Littlewood}\atop\mbox{$t$}$}}
%
\put(48,55){\vector(-2,-1){29}}
\put(72,55){\vector(2,-1){29}}
\put(20,30){\vector(2,-1){30}}
\put(100,30){\vector(-2,-1){30}}
\put(60,50){\vector(0,-1){30}}
\put(15,47.5){\makebox(20,7.5){$t=q^{\beta}$, $q\rightarrow 1$}}
\put(80,47.5){\makebox(20,7.5){$q\rightarrow 0$}}
\put(17.5,16){\makebox(20,7.5){$\beta\rightarrow 1$}}
\put(80,16){\makebox(20,7.5){$t\rightarrow 0$}}
\put(56,30){\makebox(20,7.5){$q=t$}}
\end{picture}
}

\section{$W_N$ Algebra}

In this section we recapitulate the $W_N$ algebra\cite{rFL} from the free
filed realization point of view. For comprehensive review of $W$ algebras,
see \cite{rBS}.

The $W_N$ algebra ($A_{N-1}$ type $W$ algebra) is an associative
algebra generated by spin $k$ currents\footnote{
Bar ( $\bar{}$ ) indicates non-$q$-deformed quantity.}
$\displaystyle \bar{W}^k(z)=\sum_{n\in\bz}\bar{W}^k_nz^{-n-k}$ 
($k=2,\cdots,N$) and the central charge $c$.
To define its relations it is convenient to use a free field realization.
Let us introduce free bosons 
$\displaystyle \bar{h}^i(z)=Q_h^i+h^i_0\log z-\sum_{n\neq 0}
\sfrac{1}{n}\bar{h}^i_nz^{-n}$ ($i=1,\cdots,N$), 
whose relations are given by
\ba
  \lbrack\bar{h}^i_n,\bar{h}^j_m\rbrack
  &\!\!=\!\!&
  (\delta^{ij}-\sfrac{1}{N})n\delta_{n+m,0},\quad
  \bar{h}^i_0=h^i_0, \n
  \lbrack\bar{h}^i_n,Q_h^j\rbrack
  &\!\!=\!\!&
  (\delta^{ij}-\sfrac{1}{N})\delta_{n0},\quad
  \lbrack Q_h^i,Q_h^j\rbrack=0,\quad
  \sum_{i=1}^N\bar{h}^i_n=0,\quad \sum_{i=1}^NQ_h^i=0.
\ea
These $\bar{h}^i(z)$'s correspond to the weights of the vector representation
of $A_{N-1}$ algebra, $h_i$.
If we introduce an orthonormal set of free bosons 
$\varphi(z){}^t\varphi(w)\sim{\bf 1}\log(z-w)$, and
project them on the root space, 
$\phi(z)=\varphi(z)\Bigl|_{\mbox{root space}}$,
then $\bar{h}^i$ is expressed as $\bar{h}^i(z)=h_i\cdot\phi(z)$.
We will also use free bosons $\phi^a(z)=\alpha^a\cdot\phi(z)$, 
which correspond to the simple roots $\alpha^a$ ($a=1,\cdots,N-1$).
$\phi^a(z)$ has the form 
$\displaystyle \phi^a(z)=Q_{\alpha}^a+\alpha^a_0\log z-\sum_{n\neq 0}
 \sfrac{1}{n}\bar{\alpha}^a_nz^{-n}$ 
and $\alpha^a=h_a-h_{a+1}$ implies
\be
  \lbrack\bar{\alpha}^a_n,\bar{\alpha}^b_m\rbrack=A_{ab}n\delta_{n+m,0},\quad
  \bar{\alpha}^a_0=\alpha^a_0,\quad
  \lbrack\bar{\alpha}^a_n,Q_{\alpha}^b\rbrack=A_{ab}\delta_{n0},\quad
  \lbrack Q_{\alpha}^a,Q_{\alpha}^b\rbrack=0,
\ee
where $A=(A_{ab})$ is the Cartan matrix of $A_{N-1}$ algebra.
Note that $h_i\cdot h_j=\delta_{ij}-\frac{1}{N}$,
$\alpha^a\cdot\alpha^b=A_{ab}$ and 
$\Lambda_a\cdot\Lambda_b=(A^{-1})_{ab}$, 
where $\Lambda_a=\sum_{i=1}^ah_i$ is the fundamental weight, 
$\alpha^a\cdot\Lambda_b=\delta_{ab}$.

The following Miura transformation gives the realization of the $W_N$ 
algebra (this is just the definition of the $W_N$ 
algebra by ref.\cite{rFL});
\be
  :\!\Bigl(\alpha_0\partial_z+\partial\bar{h}^1(z)\Bigr)
   \Bigl(\alpha_0\partial_z+\partial\bar{h}^2(z)\Bigr)\cdots
   \Bigl(\alpha_0\partial_z+\partial\bar{h}^N(z)\Bigr)\!:\;
  =
  \sum_{k=0}^N\bar{W}^k(z)\Bigl(\alpha_0\partial_z\Bigr)^{N-k}.
  \label{Miura}
\ee
Here $:*:$ stands for the normal ordering (i.e., non-negative mode oscillators
are moved to the right for negative mode oscillators and $Q$), and 
$\alpha_0$ is a parameter which determines the central charge.
For later convenience we parameterize $\alpha_0$ as
\be
  \alpha_0=\sqrt{\beta}-\frac{1}{\sqrt{\beta}}.
\ee
We have $W^0(z)=1$, $W^1(z)=0$, and 
\be
  -W^2(z)=L(z)=\sum_{n\in\bz}L_nz^{-n-2}=
  \sfrac12:\partial\phi(z)\cdot\partial\phi(z)\!:
  +\alpha_0\rho\cdot\partial^2\phi(z).
  \label{Vir}
\ee
This $L(z)$ generates the Virasoro algebra
\be
  [L_n,L_m]=(n-m)L_{n+m}+\frac{c}{12}(n^3-n)\delta_{n+m,0},
\ee
 with the central charge,
\be
  c=N-1-12\alpha_0^2\rho^2,
  \label{c}
\ee
where $\rho$ is the half sum of the positive roots 
$\rho=\sum_{a=1}^{N-1}\Lambda_a$ and $\rho^2=\sfrac{1}{12}N(N^2-1)$.

$W_N$ generators satisfy quadratic relations. Symbolically it is
\be
  [\bar{W}^k_n,\bar{W}^\ell_m]=\sum((\bar{W}\bar{W})+\bar{W}+1).
  \label{WW}
\ee
Here quadratic terms should be normal-ordered; For operators 
$\displaystyle A(z)=\sum_{n\in\bz-h_A}A_nz^{-n-h_A}$ and
$\displaystyle B(z)=\sum_{n\in\bz-h_B}B_nz^{-n-h_B}$, 
normal ordering $(AB)(z)$ is defined by
\ba
  (AB)(w)
  &\!\!=\!\!&
  \oint_w\frac{dz}{2\pi i}\frac{1}{z-w}A(z)B(w) \n
  &\!\!=\!\!&
  \oint_0\frac{dz}{2\pi iz}
  \Biggl(\frac{1}{1-\frac{w}{z}}A(z)B(w)
  +\frac{\frac{z}{w}}{1-\frac{z}{w}}B(w)A(z)\Biggr)
  \label{(AB)} \\
  &\!\!=\!\!&
  \sum_{n\in\bz-h_A-h_B}\biggl(
  \sum_{m\leq -h_A}A_mB_{n-m}+\sum_{m>-h_A}B_{n-m}A_m\biggr)
  \cdot w^{-n-h_A-h_B}. \nonumber
\ea

We consider the highest weight representations. The highest weight
state $|\mbox{hws}\bar{\rangle}$ is characterized by
\be
  \bar{W}^k_n|\mbox{hws}\bar{\rangle}=0\; (n>0),\quad
  \bar{W}^k_0|\mbox{hws}\bar{\rangle}=\bar{w}_k|\mbox{hws}\bar{\rangle}\; 
  (\bar{w}_k\in\bC).
\ee
The Verma module is obtained by successive action of $\bar{W}^k_{-n}$ ($n>0$) 
on $|\mbox{hws}\bar{\rangle}$. We remark that only finite number of terms
in $(WW)$ of \eq{WW} survives on $|\mbox{hws}\bar{\rangle}$.
If there exists a singular vector $|\bar{\chi}\bar{\rangle}$, which is
characterized by
\be
  \bar{W}^k_n|\bar{\chi}\bar{\rangle}=0\; (n>0),\quad
  \bar{W}^k_0|\bar{\chi}\bar{\rangle}=
  (\bar{w}_k+\bar{N}_k)|\bar{\chi}\bar{\rangle}=0\; 
  (\bar{N}_k\in\bC),
\ee
then the Verma module is reducible. To obtain an irreducible module,
the submodule on $|\bar{\chi}\bar{\rangle}$ has to be factored out.

In the following we consider the representations realized in the boson
Fock space. The Fock vacuum $|\gamma\bar{\rangle}$ is characterized by
\be
  \bar{\alpha}^a_n|\gamma\bar{\rangle}=0\; (n>0),\quad
  \alpha^a_0|\gamma\bar{\rangle}=\gamma^a|\gamma\bar{\rangle},
\ee
where $\gamma=\sum_{a=1}^{N-1}\gamma^a\Lambda_a$ ($\gamma^a\in\bC$).
$|\gamma\bar{\rangle}$ can be obtained from $|0\bar{\rangle}$ 
($\bar{\alpha}_n|0\bar{\rangle}=0$ for $n\geq 0$),
\be
  |\gamma\bar{\rangle}=
  \exp\Bigl(\sum_{a=1}^{N-1}\gamma^aQ_{\Lambda}^a\Bigr)\cdot
  |0\bar{\rangle},\quad
  Q_{\Lambda}^a=\sum_{j=1}^aQ_h^j.
\ee
The Fock space $\bar{F}_{\gamma}$ is a linear span of 
$\bar{\alpha}^{a_1}_{-n_1}\bar{\alpha}^{a_2}_{-n_2}\cdots|\gamma\bar{\rangle}$ 
($n_1,n_2,\cdots>0$).
Dual Fock space is defined similarly ( 
$\bar{\langle}\gamma|\bar{\alpha}^a_n=\bar{\langle}\gamma|\gamma^a\delta_{n0}$ 
($n\leq 0$), with normalization $\bar{\langle}\gamma|\gamma'\bar{\rangle}=
\delta_{\gamma\gamma'}$).
For studying representations in the Fock space, the most important 
tool is the screening current.
Let us introduce the screening currents $\bar{S}^a_{\pm}(z)$ 
($a=1,\cdots N-1$),
\be
  \bar{S}^a_{\pm}(z)=\;:\!e^{\alpha_{\pm}\phi^a(z)}\!:,\qquad
  \alpha_+=\sqrt{\beta},\;\;\alpha_-=\frac{-1}{\sqrt{\beta}}.
\ee
The operator product expansion between $W_N$ currents and the screening
currents is
\be
  :\!\Bigl(\alpha_0\partial_z+\partial\bar{h}^1(z)\Bigr)\cdots
   \Bigl(\alpha_0\partial_z+\partial\bar{h}^N(z)\Bigr)\!:
  \bar{S}^a_{\pm}(w)
  \sim
  \frac{\partial}{\partial w}\Bigl(\cdots\Bigr)+\mbox{reg.}
\ee
Therefore the screening charge $\oint dz\bar{S}^a_{\pm}(z)$ commutes
with $W_N$ algebra,
\be
  \Bigl[\:W_N,\oint dz\bar{S}^a_{\pm}(z)\:\Bigr]=0.
\ee
We must comment that this proposition holds only on 
a suitable state on which the contour closes.
To study representation in more detail, we have to construct the BRST charge
and analyze the BRST cohomology\cite{rFe}. Here we do not get into
this direction any more.

The Fock vacuum $|\gamma\bar{\rangle}$ satisfies the highest weight state
condition of the $W_N$ algebra. If we choose $\gamma$ as special values,
the Verma module contains singular vectors.
Let us define $\alpha_{rs}^{\pm}$ and $\tilde{\alpha}_{rs}^{\pm}$ as
\ba
  \alpha_{rs}^+
  &\!\!=\!\!&
  \sum_{a=1}^{N-1}\Bigl(\sqrt{\beta}(1+r_a-r_{a-1})
  +\frac{-1}{\sqrt{\beta}}(1+s_a)\Bigr)\Lambda_a,\quad
  \tilde{\alpha}_{rs}^+=\alpha_{rs}^+
  -\sqrt{\beta}\sum_{a=1}^{N-1}r_a\alpha^a, \n
  \alpha_{rs}^-
  &\!\!=\!\!&
  \sum_{a=1}^{N-1}\Bigl(\frac{-1}{\sqrt{\beta}}(1+r_a-r_{a-1})
  +\sqrt{\beta}(1+s_a)\Bigr)\Lambda_a,\quad
  \tilde{\alpha}_{rs}^-=\alpha_{rs}^-
  -\frac{-1}{\sqrt{\beta}}\sum_{a=1}^{N-1}r_a\alpha^a,
  \label{alphars}
\ea
where $r_a$, $s_a$ are positive integers 
such that  $r_1>r_2>\cdots>r_{N-1}$ and $r_0=0$.
The Verma module on $|\alpha_{rs}^+\rangle$ contains a singular vector 
$|\bar{\chi}_{rs}^+\bar{\rangle}$,

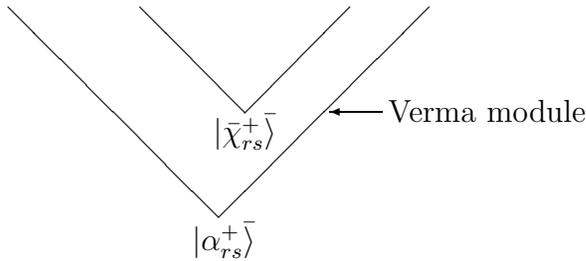
\setlength{\unitlength}{0.7mm}
\begin{picture}(140,60)
\put(50,10){\line(1,1){40}}
\put(50,10){\line(-1,1){40}}
\put(55,30){\line(1,1){20}}
\put(55,30){\line(-1,1){20}}
\put(81,30){\vector(-1,0){10}}
\put(41,2){\makebox(20,7.5){$|\alpha^+_{rs}\bar{\rangle}$}}
\put(45,22.5){\makebox(20,7.5){$|\bar{\chi}^+_{rs}\bar{\rangle}$}}
\put(91,26.5){\makebox(20,7.5){Verma module}}
\end{picture}
\ba
  |\bar{\chi}^+_{rs}\bar{\rangle}
  &\!\!=\!\!&
  \oint\prod_{a=1}^{N-1}\prod_{j=1}^{r_a}\frac{dx^a_j}{2\pi i}\cdot
  \bar{S}^1_+(x^1_1)\cdots\bar{S}^1_+(x^1_{r_1})\cdots
  \bar{S}^{N-1}_+(x^{N-1}_1)\cdots\bar{S}^{N-1}_+(x^{N-1}_{r_{N-1}})
  |\tilde{\alpha}^+_{rs}\bar{\rangle} \n
  &\!\!=\!\!&
  \oint\prod_{a=1}^{N-1}\prod_{j=1}^{r_a}\frac{dx^a_j}{2\pi ix^a_j}\cdot
  \prod_{a=1}^{N-1}\bar{\pi}\Bigl(\frac{1}{x^a},x^{a+1}\Bigr)\bar{\Delta}(x^a)
  \prod_{j=1}^{r_a}\Bigl(x^a_j\Bigr)^{-s_a}
  \Bigl[\bar{S}^a_+(x^a_j)\Bigr]_-\cdot
  |\alpha^+_{rs}\bar{\rangle},
\ea
where $\bar{\Delta}(x)$ and $\bar{\pi}(x,y)$ are
\be
  \bar{\Delta}(x)=\prod_{i\neq j}\Bigl(1-\frac{x_i}{x_j}\Bigr)^{\beta},\quad
  \bar{\pi}(x,y)=\prod_{i,j}(1-x_iy_j)^{-\beta},
\ee
and $[*]_-$ stands for the negative mode oscillator part.

This singular vector is related to the Jack symmetric polynomial.
To see this, we consider a map from the Fock space to the ring of 
symmetric function;
\be
  \bar{F}_{\gamma}\ni|f\bar{\rangle}\mapsto 
  f(x)=\bar{\langle}\gamma|
  \exp\Bigl(\sqrt{\beta}\sum_{n>0}\sfrac{1}{n}\bar{h}^1_np_n\Bigr)
  |f\bar{\rangle},
\ee
where $p_n$ is a power sum symmetric polynomial $p_n=\sum_i(x_i)^n$.
By this map, $\bar{\alpha}^a_{-n}$ is replaced by 
$\delta^{a1}\sqrt{\beta}p_n$.
We remark that $\exp(\sqrt{\beta}\sum_{n>0}\sfrac{1}{n}\bar{h}^1_nz^{-n})$ 
is the positive oscillator part of the vertex operator corresponding to
the vector representation.
Then the singular vector $|\bar{\chi}^+_{rs}\bar{\rangle}$ is mapped to the
Jack symmetric polynomial,
\ba
  &&
  \bar{\langle}\alpha^+_{rs}|\exp\Bigl(
  \sqrt{\beta}\sum_{n>0}\sfrac{1}{n}\bar{h}^1_np_n\Bigr)
  |\bar{\chi}^+_{rs}\bar{\rangle} \n
  &=\!\!&
  \oint\prod_{a=1}^{N-1}\prod_{j=1}^{r_a}\frac{dx^a_j}{2\pi ix^a_j}\cdot
  \bar{\pi}(x,x^1)
  \prod_{a=1}^{N-1}\bar{\pi}\Bigl(\frac{1}{x^a},x^{a+1}\Bigr)\bar{\Delta}(x^a)
  \prod_{j=1}^{r_a}\Bigl(x^a_j\Bigr)^{-s_a} \\
  &\!\!\propto\!\!&
  J_{\lambda}(x;\beta). \nonumber
\ea
In the last equation we have changed the integration variable 
$x^a\rightarrow\frac{1}{x^a}$ and used the integral representation
of the Jack symmetric polynomial\cite{rAMOS}.
Here the partition $\lambda$ is 
$\lambda'=((r_1)^{s_1},(r_2)^{s_2},\cdots,(r_{N-1})^{s_{N-1}})$,
namely, corresponds to the following Young diagram
\def\generalYoung{
\vskip.25cm
\noindent
\makebox[  4cm]{ }
\makebox[  2cm]{$s_1$}\hskip-.4pt
\makebox[1.7cm]{$s_2$}
\makebox[1.4cm]{$\!\cdots\cdot$}
\makebox[1.4cm]{$\!\!s_{N-2}$}\hskip-.5pt
\makebox[1.3cm]{$\!\!s_{N-1}$}
\hfill\break
 \makebox[  4cm][r]{$\hfill\lambda\:=\;$}
\framebox[  2cm][l]{\rule[  -1cm]{0cm}{  2cm}
                    \raisebox{.2cm}{$\!r_1$}}\hskip-.4pt
\framebox[1.7cm][l]{\rule[-0.7cm]{0cm}{1.7cm}
                    \raisebox{.2cm}{$\!r_2$}}
\hskip-0.15cm\rule[1.105cm]{1.6cm}{0.4pt}\hskip-1.55cm
 \makebox[1.4cm]   {\raisebox{.2cm}{$\,\,\cdots\cdot$}}
\framebox[1.4cm][l]{\rule[-0.3cm]{0cm}{1.3cm}
                    \raisebox{.2cm}{$\!r_{N-2}$}}\hskip-.4pt
\framebox[1.3cm][l]{\rule[0.0cm]{0cm}{1.0cm}
                    \raisebox{.2cm}{$\!r_{N-1}$}}
\makebox[0.2cm][r]{.}
\vskip.4cm
}
\generalYoung

\noindent
Therefore, in the free field realization, the singular vector
of the $W_N$ algebra realizes the Jack symmetric polynomial
with the Young diagram composed of $N-1$ rectangles.
Similarly we have singular vectors 
$|\bar{\chi}^-_{rs}\bar{\rangle}=
 \oint \bar{S}_-\cdots\bar{S}_-|\tilde{\alpha}^-_{rs}\bar{\rangle}$,
and another type of integral representation of the Jack polynomial\cite{rAMOS}.

\section{Quantum Deformed $W_N$ Algebra}

In this section we define and explain the quantum deformed $W_N$ 
algebra\cite{rAKOS,rFF} along the line of $W_N$ in section 2.

The quantum deformed $W_N$ algebra, $q$-$W_N$, is an associative algebra
generated by 
$\displaystyle W^i(z)=\sum_{n\in\bz}W^i_nz^{-n}$ ($i=1,\cdots,N-1$) 
and contains two parameters $q$ and $t$ ($q,t\in\bC$). 
We will often use the following notation,
\be
  p=qt^{-1},\quad 
  q=e^{\hbar}=e^{\frac{1}{\sqrt{\beta}}\hbar'},\quad
  t=q^{\beta}=e^{\sqrt{\beta}\hbar'}.
\ee
To define the relations of $q$-$W_N$ generators, we first introduce
fundamental bosons $h^i_n$ ($n\in\bZ$) and $Q^i_h$ ($i=1,\cdots,N$),
\ba
  \lbrack h^i_n,h^j_m\rbrack
  &\!\!=\!\!&
  \frac{1}{n}
  (q^{\frac{n}{2}}-q^{-\frac{n}{2}})(t^{\frac{n}{2}}-t^{-\frac{n}{2}})
  \frac{p^{ \frac{n}{2}N(\delta^{ij}-\frac{1}{N})}-
        p^{-\frac{n}{2}N(\delta^{ij}-\frac{1}{N})}}
       {p^{\frac{n}{2}N}-p^{-\frac{n}{2}N}}
  p^{\frac{n}{2}N\mbox{sgn}(j-i)}\delta_{n+m,0}, \n
  \lbrack h^i_n,Q_h^j\rbrack
  &\!\!=\!\!&
  (\delta^{ij}-\sfrac{1}{N})\delta_{n0},\quad
  \lbrack Q_h^i,Q_h^j\rbrack=0,\quad
  \sum_{i=1}^Np^{in}h^i_n=0, \quad \sum_{i=1}^NQ_h^i=0,
  \label{hh}
\ea
where $\mbox{sgn}(x)$ is 
$\mbox{sgn}(x)=1(\mbox{for }x>0),0(\mbox{for }x=0),-1(\mbox{for }x<0)$.
We remark that the fractional part containing $p$ is a $p$-analogue
of $\delta^{ij}-\frac{1}{N}$.
Let us define transformations $\theta$, $\omega$, and 
$\omega'=\theta\omega$ as follows:
\be
\begin{tabular}{|c||c|c|c|}
  \hline
  $\qquad\qquad$&$\qquad\theta\qquad$&$\qquad\omega\qquad$&
  $\qquad\omega'\qquad$ \\ \hline
  $q$&$q^{-1}$&$t$&$t^{-1}$ \\
  $t$&$t^{-1}$&$q$&$q^{-1}$ \\
  $h^i_n$ ($n\neq 0$)&$h^{N+1-i}_n$&$h^{N+1-i}_n$&$h^i_n$ \\
  $h^i_0$&$-h^{N+1-i}_0$&$h^{N+1-i}_0$&$-h^i_0$ \\
  $Q_h^i$&$-Q_h^{N+1-i}$&$Q_h^{N+1-i}$&$-Q_h^i$ \\
  \hline
\end{tabular}
\ee
For example $\theta\cdot q=q^{-1}$, $\omega\cdot h^i_0=h^{N+1-i}_0$.
These are involutions, $\theta^2=\omega^2=\omega^{\prime\,2}=1$.
Then \eq{hh} is invariant under $\theta$, $\omega$ and $\omega'$.
For later convenience we list up transformation rules for various
quantities (for their definitions see below):
\be
\begin{tabular}{|c||c|c|c|}
  \hline
  $\qquad\qquad$&$\qquad\theta\qquad$&$\qquad\omega\qquad$&
  $\qquad\omega'\qquad$ \\ \hline
  $\hbar'$&$-\hbar'$&$\hbar'$&$-\hbar'$ \\
  $\beta$&$\beta$&$\beta^{-1}$&$\beta^{-1}$ \\
  $p$&$p^{-1}$&$p^{-1}$&$p$ \\
  $\alpha^a_n$ ($n\neq 0$)&$-\alpha^{N-a}_n$&$-\alpha^{N-a}_n$&$\alpha^a_n$ \\
  $\alpha^a_0$&$\alpha^{N-a}_0$&$-\alpha^{N-a}_0$&$-\alpha^a_0$ \\
  $Q_{\alpha}^a$&$Q_{\alpha}^{N-a}$&$-Q_{\alpha}^{N-a}$&$-Q_{\alpha}^a$ \\
  $\Lambda_i(z)$&$\Lambda_{N+1-i}(z)$&$\Lambda_{N+1-i}(z)$&$\Lambda_i(z)$ \\
  $W^i(z)$&$W^i(z)$&$W^i(z)$&$W^i(z)$ \\
  $f^{ij}(x)$&$f^{ij}(x)$&$f^{ij}(x)$&$f^{ij}(x)$ \\
  $S_{\pm}^a(z)$&$S_{\pm}^{N-a}(z)$&$S_{\mp}^{N-a}(z)$&$S_{\mp}^a(z)$ \\
  \hline
\end{tabular}
\ee

In the case of $W_N$ algebra, the building block of $W_N$ currents
is a boson $\bar{h}^i(z)$, and the Miura transformation is an equation
of differential operator.
On the other hand the building block of $q$-$W_N$ currents is an exponentiated
boson $\Lambda_i(z)$ ($i=1,\cdots,N$),
\be
  \Lambda_i(z)=\;:\exp\Bigl(\sum_{n\neq 0}h^i_nz^{-n}\Bigr)\!:
  q^{\sqrt{\beta}h^i_0}p^{\frac{N+1}{2}-i},
\ee
and the differential operator (shift operator) $\partial_z$  is replaced 
by $p$-difference operator ($p$-shift operator) $p^{D_z}$ 
($D_z=z\partial_z$); $p^{D_z}f(z)=f(pz)$.
$q$-deformed Miura transformation is given by
\be
  :\!\Bigl(p^{D_z}-\Lambda_1(z)\Bigr)
   \Bigl(p^{D_z}-\Lambda_2(p^{-1}z)\Bigr)\cdots
   \Bigl(p^{D_z}-\Lambda_N(p^{1-N}z)\Bigr)\!:\;
  =
  \sum_{i=0}^N(-1)^iW^i(p^{\frac{1-i}{2}}z)p^{(N-i)D_z}.
  \label{q-Miura}
\ee
{}From this, $W^i(z)$ is expressed as
\be
  W^i(z)=\sum_{1\leq j_1<j_2<\cdots<j_i\leq N}
  :\!\Lambda_{j_1}(p^{\frac{i-1}{2}}z)
   \Lambda_{j_2}(p^{\frac{i-3}{2}}z)\cdots
   \Lambda_{j_i}(p^{-\frac{i-1}{2}}z)\!:,
\ee
and $W^0(z)=W^N(z)=1$.
$W^i(z)$ is invariant under $\theta$, $\omega$ and $\omega'$.
Eq.\eq{q-Miura} is equivalent to 
$$
  :\!\Bigl(p^{-D_z}-\Lambda_N(z)\Bigr)
   \Bigl(p^{-D_z}-\Lambda_{N-1}(pz)\Bigr)\cdots
   \Bigl(p^{-D_z}-\Lambda_1(p^{N-1}z)\Bigr)\!:\;
  =
  \sum_{i=0}^N(-1)^iW^i(p^{\frac{i-1}{2}}z)p^{-(N-i)D_z},
$$
\vskip -5mm
\ba
  \!\!\!\!\!\!\!\!
  :\!\Bigl(1-\Lambda_1(z)p^{-D_z}\Bigr)
   \Bigl(1-\Lambda_2(z)p^{-D_z}\Bigr)\cdots
   \Bigl(1-\Lambda_N(z)p^{-D_z}\Bigr)\!:
  &\!\!\!=\!\!\!&
  \sum_{i=0}^N(-1)^iW^i(p^{\frac{1-i}{2}}z)p^{-iD_z}\!\!,\;\;\\
  :\!\Bigl(1-\Lambda_N(z)p^{D_z}\Bigr)
   \Bigl(1-\Lambda_{N-1}(z)p^{D_z}\Bigr)\cdots
   \Bigl(1-\Lambda_1(z)p^{D_z}\Bigr)\!:
  &\!\!\!=\!\!\!&
  \sum_{i=0}^N(-1)^iW^i(p^{\frac{i-1}{2}}z)p^{iD_z}.\nonumber
\ea

The relation of $q$-$W_N$ algebra is quadratic. Symbolically it is
\be
  f^{ij}(\sfrac{w}{z})W^i(z)W^j(w)-W^j(w)W^i(z)f^{ji}(\sfrac{z}{w})
  =
  \sum\Bigl(\no WW\no+W+1\Bigr).
  \label{fWW}
\ee
We remark that this relation is invariant under $\theta$, $\omega$ 
and $\omega'$.
The structure function 
$\displaystyle f^{ij}(x)=\sum_{\ell=0}^{\infty}f^{ij}_{\ell}x^{\ell}$ 
is defined by
\ba
  \!\!\!\!\!
  &&f^{ij}(x)=
  \exp\biggl(-\sum_{n>0}\frac{1}{n}
  (q^{\frac{n}{2}}-q^{-\frac{n}{2}})(t^{\frac{n}{2}}-t^{-\frac{n}{2}}) \n
  \!\!\!\!\!
  &&\qquad\qquad\qquad\quad\quad\times
  \frac{p^{\frac{n}{2}\min(i,j)}-p^{-\frac{n}{2}\min(i,j)}}
       {p^{\frac{n}{2}}-p^{-\frac{n}{2}}}
  \frac{p^{\frac{n}{2}(N-\max(i,j))}-p^{-\frac{n}{2}(N-\max(i,j))}}
       {p^{\frac{n}{2}N}-p^{-\frac{n}{2}N}}
  x^n\biggl).
\ea
Note that the fractional part containing $p$ is a $p$-analogue of
$(A^{-1})_{ij}$ and $f^{ij}(x)=f^{ji}(x)=f^{N-i,N-j}(x)$.
Quadratic terms have to be normal-ordered,
\ba
  &\!\!&\!\!\!\!
  \no W^i(rw)W^j(w)\no \n
  &\!\!=\!\!&
  \oint\frac{dz}{2\pi iz}\Biggl(
  \frac{1}{1-\frac{rw}{z}}f^{ij}(\sfrac{w}{z})W^i(z)W^j(w)
  +\frac{\frac{z}{rw}}{1-\frac{z}{rw}}W^j(w)W^i(z)f^{ji}(\sfrac{z}{w})
  \Biggr) \n
  &\!\!=\!\!&
  \sum_{n\in\bz}\sum_{m=0}^{\infty}\sum_{\ell=0}^m f^{ij}_\ell \Bigl(
  r^{m-\ell}W^i_{-m}W^j_{n+m}+r^{\ell-m-1}W^j_{n-m-1}W^i_{m+1}
  \Bigr)\cdot w^{-n},
\ea
where $(1-x)^{-1}$ stands for $\displaystyle \sum_{n\geq 0}x^n$.
This normal ordering $\no*\no$ is a generalization of \eq{(AB)}.

Here we present some examples of them.
The relation of $W^1(z)$ and $W^j(z)$ for $j\geq 1$ is
\ba
  &\!\!\!\!\!\!&\!\!\!\!\!\!
  f^{1j}(\sfrac{w}{z})W^1(z)W^j(w) - W^j(w)W^1(z)f^{j1}(\sfrac{z}{w}) \n
  &\!\!=\!\!&
  -\frac{(1-q)(1-t^{-1})}{1-p}\biggl(
  \delta(p^{\frac{j+1}{2}}\sfrac{w}{z})W^{j+1}(p^{\frac12}w)
  -\delta(p^{-\frac{j+1}{2}}\sfrac{w}{z})W^{j+1}(p^{-\frac12}w)
  \biggr),
\ea
with $\displaystyle \delta(x)=\sum_{n\in\bZ} x^n$;
and that of $W^2(z)$ and $W^j(z)$ for $j\geq 2$ is
\ba
  &\!\!\!\!\!\!&\!\!\!\!\!\!
  f^{2j}(\sfrac{w}{z})W^2(z)W^j(w)-W^j(w)W^2(z)f^{j2}(\sfrac{z}{w}) \n
  &\!\!=\!\!&
  -\frac{(1-q)(1-t^{-1})}{1-p}\frac{(1-qp)(1-t^{-1}p)}{(1-p)(1-p^2)}
  \biggl(\delta(p^{\frac{j}{2}+1}\sfrac{w}{z})W^{j+2}(pw) 
  -\delta(p^{-\frac{j}{2}-1}\sfrac{w}{z})W^{j+2}(p^{-1}w)\biggr)\n
  &\!\!\!\!\!\!&
  -\frac{(1-q)(1-t^{-1})}{1-p}
  \biggl(\delta(p^{\frac{j}{2}}\sfrac{w}{z})
  \no W^1(p^{-\frac{1}{2}}z)W^{j+1}(p^{\frac{1}{2}}w)\no 
  -\delta(p^{-\frac{j}{2}}\sfrac{w}{z})
  \no W^1(p^{\frac{1}{2}}z)W^{j+1}(p^{-\frac{1}{2}}w)\no \biggr) \n
  &\!\!\!\!\!\!&
  +\frac{(1-q)^2(1-t^{-1})^2}{(1-p)^2}
  \biggl(\delta(p^{\frac{j}{2}}\sfrac{w}{z})
  \Bigl(\frac{p^2}{1-p^2}W^{j+2}(pw)+\frac{1}{1-p^j}W^{j+2}(w)\Bigr) \\
  &\!\!\!\!\!\!&\hspace{35mm}
  -\delta(p^{-\frac{j}{2}}\sfrac{w}{z})
  \Bigl(\frac{p^j}{1-p^j}W^{j+2}(w)+\frac{1}{1-p^2}W^{j+2}(p^{-1}w)\Bigr)
  \biggr), \nonumber
\ea
with $W^i(z)=0$ for $i>N$.
The main term of
$$
  f^{ij}(\sfrac{w}{z})W^i(z)W^j(w)-W^j(w)W^i(z)f^{ji}(\sfrac{z}{w})\qquad
  (i\leq j)
$$
is
\ba
  &\!\!\!\!\!\!&
  -\frac{(1-q)(1-t^{-1})}{1-p}\sum_{k=1}^{\min(i,N-j)}\prod_{\ell=1}^{k-1}
   \frac{(1-qp^\ell)(1-t^{-1}p^\ell)}{(1-p^\ell)(1-p^{\ell+1})} \\
  &\!\!\!\!\!\!&\!\!\!\!\!\! \times
  \biggl(\delta(p^{\frac{j-i}{2}+k}\sfrac{w}{z})
  \no W^{i-k}(p^{-\frac{k}{2}}z)W^{j+k}(p^{\frac{k}{2}}w)\no
  -\delta(p^{-\frac{j-i}{2}-k}\sfrac{w}{z})
  \no W^{i-k}(p^{\frac{k}{2}}z)W^{j+k}(p^{-\frac{k}{2}}w)\no
  \biggr). \nonumber
\ea
By taking $\beta\rightarrow 0$ limit ($q$ fixed), this main term 
reduces to $q$-deformed $W_N$ Poisson bracket algebra in \cite{rFR}.

The highest weight state $|\mbox{hws}\rangle$ is characterized by
\be
  W^i_n|\mbox{hws}\rangle=0\; (n>0),\quad
  W^i_0|\mbox{hws}\rangle=w_i|\mbox{hws}\rangle\; (w_i\in\bC).
\ee
We remark that only finite number of terms in $\no WW\no$ of \eq{fWW} 
are non-vanishing on $|\mbox{hws}\rangle$.
The boson Fock space $F_{\gamma}$ is defined like as in section 2.
To construct screening currents we define root bosons $\alpha^a_n$ 
($n\in\bZ$) and $Q_{\alpha}^a$ ($a=1,\cdots,N-1$) as
$\alpha^a_n=h^a_n-h^{a+1}_n$, 
$Q_{\alpha}^a=Q_h^a-Q_h^{a+1}$.
They satisfy
\ba
  \lbrack \alpha^a_n,\alpha^b_m\rbrack
  &\!\!=\!\!&
  \frac{1}{n}
  (q^{\frac{n}{2}}-q^{-\frac{n}{2}})(t^{\frac{n}{2}}-t^{-\frac{n}{2}})
  \frac{p^{\frac{n}{2}A_{ab}}-p^{-\frac{n}{2}A_{ab}}}
       {p^{\frac{n}{2}}-p^{-\frac{n}{2}}}
  p^{\frac{n}{2}\mbox{sgn}(b-a)}\delta_{n+m,0}, \n
  \lbrack \alpha^a_n,Q_{\alpha}^b\rbrack
  &\!\!=\!\!&
  A_{ab}\delta_{n0},\quad \lbrack Q_{\alpha}^a,Q_{\alpha}^b\rbrack=0.
\ea
The fractional part containing $p$ is a $p$-analogue of $A_{ab}$.
Screening currents are defined by
\ba
  S_+^a(z)
  &\!\!=\!\!&
  :\exp\biggl(-\sum_{n\neq 0}
  \frac{\alpha^a_n}{q^{\frac{n}{2}}-q^{-\frac{n}{2}}}z^{-n}\biggr)\!:
  e^{\sqrt{\beta}Q_{\alpha}^a}z^{\sqrt{\beta}\alpha^a_0}, \\
  S_-^a(z)
  &\!\!=\!\!&
  :\exp\biggl(\sum_{n\neq 0}
  \frac{\alpha^a_n}{t^{\frac{n}{2}}-t^{-\frac{n}{2}}}z^{-n}\biggr)\!:
  e^{\frac{-1}{\sqrt{\beta}}Q_{\alpha}^a}
  z^{\frac{-1}{\sqrt{\beta}}\alpha^a_0}.
\ea
Commutation relation of $W^i(z)$ and $S^a_{\pm}(w)$ can be expressed as
a total difference:
\ba
  &&
  \Bigl[\::\!\Bigl(p^{D_z}-\Lambda_1(z)\Bigr)
  \Bigl(p^{D_z}-\Lambda_2(p^{-1}z)\Bigr)\cdots
  \Bigl(p^{D_z}-\Lambda_N(p^{1-N}z)\Bigr)\!:\;,
  S_{\pm}^a(w)\:\Bigr] \n
  &=\!\!&
  (q^{\frac{n}{2}}-q^{-\frac{n}{2}})(t^{\frac{n}{2}}-t^{-\frac{n}{2}})
  \frac{d}{d_{q\atop t}w}\Biggl(
  :\!\Bigl(p^{D_z}-\Lambda_1(z)\Bigr)\cdots
  \Bigl(p^{D_z}-\Lambda_{a-1}(p^{2-a}z)\Bigr) \\
  &&\hspace{15mm}
  \times w\delta(p^{a-1}\sfrac{w}{z})A_{\pm}^a(w)p^{D_z}\times
  \Bigl(p^{D_z}-\Lambda_{a+2}(p^{-1-a}z)\Bigr)
  \Bigl(p^{D_z}-\Lambda_N(p^{1-N}z)\Bigr)\!:\Biggr), \nonumber
\ea
where $A_{\pm}^a(w)$ and $\frac{d}{d_{\xi}w}$ are given by\footnote{
We have used the slightly different notation from that of \cite{rAKOS}
(say `old');\\
$S_+^a(w)=S_+^{a,\mbox{old}}(q^{-\frac12}w)q^{\frac12\sqrt{\beta}\alpha^a_0}$,
$S_-^a(w)=S_-^{a,\mbox{old}}(t^{-\frac12}w)q^{-\frac12\sqrt{\beta}\alpha^a_0}$,
$A_+^a(w)=A_+^{a,\mbox{old}}(w)q^{\frac12\sqrt{\beta}\alpha^a_0}p^{\frac12}$,
$A_-^a(w)=A_-^{a,\mbox{old}}(w)q^{-\frac12\sqrt{\beta}\alpha^a_0}p^{-\frac12}$,
$\frac{d}{d_{\xi}w}f(w)=
 \xi^{\frac12}\Bigl(\frac{d}{d_{\xi}w}\Bigr)^{\mbox{old}}f(\xi^{-\frac12}w)$
}
\ba
  A_+^a(w)
  &\!\!=\!\!&
  :\exp\biggl(\sum_{n\neq 0}
  \frac{q^{\frac{n}{2}}h^{a+1}_n-q^{-\frac{n}{2}}h^a_n}
  {q^{\frac{n}{2}}-q^{-\frac{n}{2}}}w^{-n}\biggr)\!:
  e^{\sqrt{\beta}Q_{\alpha}^a}w^{\sqrt{\beta}\alpha^a_0}
  q^{\frac12\sqrt{\beta}(h^a_0+h^{a+1}_0)}p^{\frac{N}{2}-a}, \n
  A_-^a(w)
  &\!\!=\!\!&
  \omega'\cdot A_+^a(w), \\
  \frac{d}{d_{\xi}w}f(w)
  &\!\!=\!\!&
  \frac{f(\xi^{\frac12}w)-f(\xi^{-\frac12}w)}
  {(\xi^{\frac12}-\xi^{-\frac12})w}
  =\frac{d}{d_{\xi^{-1}}w}f(w).
\ea
Then screening charge $\oint dzS_{\pm}^a(z)$ commutes with $q$-$W_N$ algebra,
\be
  \Bigl[\:\mbox{$q$-$W_N$},\oint dzS_{\pm}^a(z)\:\Bigr]=0.
\ee
Exactly speaking we have to include the zero mode 
factor\cite{rLP2,rAJMP,rJLMP} and study BRST cohomology.

The Fock vacuum $|\gamma\rangle$ satisfies the highest weight condition
of $q$-$W_N$. In the Verma module generated by 
$|\alpha^+_{rs}\rangle$ ($\alpha^+_{rs}$ is given in \eq{alphars})
we have a singular vector $|\chi^+_{rs}\rangle$, 
\ignore{
\setlength{\unitlength}{0.7mm}
\begin{picture}(140,60)
\put(50,10){\line(1,1){40}}
\put(50,10){\line(-1,1){40}}
\put(55,30){\line(1,1){20}}
\put(55,30){\line(-1,1){20}}
\put(81,30){\vector(-1,0){10}}
\put(41,2){\makebox(20,7.5){$|\alpha^+_{rs}\rangle$}}
\put(45,22.5){\makebox(20,7.5){$|\chi^+_{rs}\rangle$}}
\put(91,26.5){\makebox(20,7.5){Verma module}}
\end{picture}
}
\ba
  \!\!
  |\chi^+_{rs}\rangle
  &\!\!=\!\!&
  \oint\prod_{a=1}^{N-1}\prod_{j=1}^{r_a}\frac{dx^a_j}{2\pi i}\cdot
  S^1_+(x^1_1)\cdots S^1_+(x^1_{r_1})\cdots
  S^{N-1}_+(x^{N-1}_1)\cdots S^{N-1}_+(x^{N-1}_{r_{N-1}})
  |\tilde{\alpha}^+_{rs}\rangle \n
  &\!\!=\!\!&
  \oint\prod_{a=1}^{N-1}\prod_{j=1}^{r_a}\frac{dx^a_j}{2\pi ix^a_j}\cdot\!
  \prod_{a=1}^{N-1}\pi\Bigl(\frac{1}{x^a},px^{a+1}\Bigr)\Delta(x^a)C(x^a)
  \prod_{j=1}^{r_a}\Bigl(x^a_j\Bigr)^{-s_a}
  \Bigl[S^a_+(x^a_j)\Bigr]_-\!\!\cdot
  |\alpha^+_{rs}\rangle,
\ea
where $\Delta(x)$, $\pi(x,y)$ and $C(x)$ are given by
\ignore{
\ba
  \!\!\!\!&&
  \Delta(x)=
  \prod_{i\neq j}\exp\biggl(-\sum_{n>0}\frac{1}{n}
  \frac{t^{\frac{n}{2}}-t^{-\frac{n}{2}}}{q^{\frac{n}{2}}-q^{-\frac{n}{2}}}
  p^{-\frac{n}{2}}\frac{x_j^n}{x_i^n}\biggr),\quad
  \pi(x,y)=
  \prod_{i,j}\exp\biggl(\sum_{n>0}\frac{1}{n}
  \frac{t^{\frac{n}{2}}-t^{-\frac{n}{2}}}{q^{\frac{n}{2}}-q^{-\frac{n}{2}}}
  p^{-\frac{n}{2}}x_i^ny_j^n\biggr), \n
  \!\!\!\!&&
  C(x)
  =
  \prod_{i<j}^r\exp\biggl(\sum_{n>0}\frac{1}{n}
  \frac{t^{\frac{n}{2}}-t^{-\frac{n}{2}}}{q^{\frac{n}{2}}-q^{-\frac{n}{2}}}
  \Bigl(p^{-\frac{n}{2}}\frac{x_i^n}{x_j^n}
        -p^{\frac{n}{2}}\frac{x_j^n}{x_i^n}\Bigr)\biggr)
  \cdot
  \prod_{i=1}^rx_i^{(r+1-2i)\beta}.
\ea
}
$$
  \Delta(x)=
  \prod_{i\neq j}\exp\biggl(-\sum_{n>0}\frac{1}{n}
  \frac{t^{\frac{n}{2}}-t^{-\frac{n}{2}}}{q^{\frac{n}{2}}-q^{-\frac{n}{2}}}
  p^{-\frac{n}{2}}\frac{x_j^n}{x_i^n}\biggr),\quad
  \pi(x,y)=
  \prod_{i,j}\exp\biggl(\sum_{n>0}\frac{1}{n}
  \frac{t^{\frac{n}{2}}-t^{-\frac{n}{2}}}{q^{\frac{n}{2}}-q^{-\frac{n}{2}}}
  p^{-\frac{n}{2}}x_i^ny_j^n\biggr),
$$
\vskip -5mm
\be
  C(x)=
  \prod_{i<j}^r\exp\biggl(\sum_{n>0}\frac{1}{n}
  \frac{t^{\frac{n}{2}}-t^{-\frac{n}{2}}}{q^{\frac{n}{2}}-q^{-\frac{n}{2}}}
  \Bigl(p^{-\frac{n}{2}}\frac{x_i^n}{x_j^n}
        -p^{\frac{n}{2}}\frac{x_j^n}{x_i^n}\Bigr)\biggr)
  \cdot
  \prod_{i=1}^rx_i^{(r+1-2i)\beta}.
\ee
We remark that $C(x)$ is a pseudo-constant under $q$-shift; 
$q^{D_{x_i}}C(x)=C(x)$ ($\forall i)$.

What we want to show is that this singular vector $|\chi^+_{rs}\rangle$ is
related to the Macdonald symmetric polynomial.
To establish this relation,
we consider a map from the Fock space to the ring of symmetric function;
\be
  F_{\gamma}\ni|f\rangle\mapsto 
  f(x)=\langle\gamma|\exp\biggl(
  \sum_{n>0}\frac{h^1_n}{q^{\frac{n}{2}}-q^{-\frac{n}{2}}}p_n\biggr)
  |f\rangle,\quad p_n=\sum_i(x_i)^n.
\ee
This map replaces $\alpha^a_{-n}$ with  
$\delta^{a1}\frac{1}{n}(t^{\frac{n}{2}}-t^{-\frac{n}{2}})p^{\frac{n}{2}}p_n$.
We remark that 
$\exp(\sum_{n>0}\frac{h^1_n}{q^{\frac{n}{2}}-q^{-\frac{n}{2}}}z^{-n})$ 
is the positive oscillator part of the vertex operator corresponding to
the vector representation.
The singular vector $|\chi^+_{rs}\rangle$ is mapped to the
Macdonald symmetric polynomial,
\ba
  &&
  \langle\alpha^+_{rs}|\exp\biggl(
  \sum_{n>0}\frac{h^1_n}{q^{\frac{n}{2}}-q^{-\frac{n}{2}}}p_n\biggr)
  |\chi^+_{rs}\rangle \n
  &=\!\!&
  \oint\prod_{a=1}^{N-1}\prod_{j=1}^{r_a}\frac{dx^a_j}{2\pi ix^a_j}\cdot
  \pi(x,px^1)
  \prod_{a=1}^{N-1}\pi\Bigl(\frac{1}{x^a},px^{a+1}\Bigr)\Delta(x^a)C(x^a)
  \prod_{j=1}^{r_a}\Bigl(x^a_j\Bigr)^{-s_a} \\
  &\propto\!\!&
  P_{\lambda}(x;q,t). \nonumber
\ea
In the last equation we have changed the integration variable 
$x^a\rightarrow(p^ax^a)^{-1}$ and used the integral representation
of the Macdonald symmetric polynomial\cite{rAMOS,rAOS}.
Here the partition $\lambda$ is same as in section 2,
$\lambda'=((r_1)^{s_1},(r_2)^{s_2},\cdots,(r_{N-1})^{s_{N-1}})$.
Therefore, in the free field realization, the singular vector
of the $q$-$W_N$ algebra realizes the Macdonald symmetric polynomial
with the Young diagram composed of $N-1$ rectangles.
Similarly we have singular vectors 
$|\chi^-_{rs}\rangle=\omega'\cdot|\chi^+_{rs}\rangle=
 \oint S_-\cdots S_-|\tilde{\alpha}^-_{rs}\rangle$, 
and another type of integral representation of the Macdonald
polynomial\cite{rAKOS}.

We have shown the condition (\romannumeral2) of $q$-deformation 
(see footnote 1). Next let us check the condition (\romannumeral1); 
the classical limit, i.e., $q\rightarrow 1$ limit 
($\hbar'\rightarrow 0$, $\beta$ fixed).
The fundamental boson $h^i_n$ can be expressed in a linear combination
of $\bar{h}^i_n$,
\be
  h^i_n=\hbar'\sum_{j=1}^Nd^{ij}_n\bar{h}^j_n=
  \hbar'\bar{h}^i_n+O(\hbar^{\prime\, 2})\quad (n\neq 0),
  \label{hin}
\ee
where $d^{ij}_n\in\bC$. Then $\Lambda_i(z)$ and $p^{D_z}$ have the
following $\hbar'$ expansion, 
\ba
  \Lambda_i(z)
  &\!\!=\!\!&
  1+\hbar'\Bigl(D\bar{h}^i(z)-(\sfrac{N+1}{2}-i)\alpha_0\Bigr)
   +O(\hbar^{\prime\, 2}), \n
  p^{D_z}
  &\!\!=\!\!&
  1-\hbar'\alpha_0D_z+O(\hbar^{\prime\, 2}), \\
  p^{D_z}-\Lambda_i(p^{1-i}z)
  &\!\!=\!\!&
  -\hbar'z^{\frac{N+1}{2}-i+1}
  \Bigl(\alpha_0\partial_z+\partial\bar{h}^i(z)\Bigr)z^{-(\frac{N+1}{2}-i)}
  +O(\hbar^{\prime\, 2}). \nonumber
\ea
Thus L.H.S. of \eq{q-Miura} is
\be
  (-\hbar')^Nz^{\frac{N+1}{2}}
  :\Bigl(\alpha_0\partial_z+\partial\bar{h}^1(z)\Bigr)
   \Bigl(\alpha_0\partial_z+\partial\bar{h}^2(z)\Bigr)\cdots
   \Bigl(\alpha_0\partial_z+\partial\bar{h}^N(z)\Bigr):\;
  z^{\frac{N-1}{2}}
  \times\Bigl(1+O(\hbar')\Bigr).
  \label{q-Miura-LHS}
\ee
This is nothing but \eq{Miura}. So $q$-$W_N$ algebra reduces to
$W_N$ algebra with central charge \eq{c}. 
However relation between $W^i(z)$ and $\bar{W}^k(z)$ is nontrivial.
R.H.S. of \eq{q-Miura} contains $\hbar^{\prime\,\ell}$ ($\ell<N$) terms
which must vanish, and $\hbar^{\prime\,N}$ term yields the $W_N$ algebra.
We will demonstrate these for explicit examples in the Appendix.

For other quantities, $q\rightarrow 1$ limit is straightforward,
\ba
  &&
  \alpha^a_n=\hbar'\bar{\alpha}^a_n\times\Bigl(1+O(\hbar)\Bigr)\quad
  (n\neq 0),\quad
  |\gamma\rangle=|\gamma\bar{\rangle},\n
  &&
  S^a_{\pm}(z)=\bar{S}^a_{\pm}(z)\times\Bigl(1+O(\hbar)\Bigr),\quad\;
  |\chi^{\pm}_{rs}\rangle=|\bar{\chi}^{\pm}_{rs}\bar{\rangle}
  \times\Bigl(1+O(\hbar)\Bigr),\\
  &&
  \Delta(x)=\bar{\Delta}(x)\times\Bigl(1+O(\hbar)\Bigr),\qquad
  \pi(x,y)=\bar{\pi}(x, y)\times\Bigl(1+O(\hbar)\Bigr), \mbox{ etc.}
  \nonumber
\ea

\section{Discussion}

We have defined a quantum deformed $W_N$ algebra, $q$-$W_N$.
There are many interesting points to be clarified in the
future study; representation theory and applications to physics.
Here we write down some of them.

\noindent(\romannumeral1)~ 
Explicit form of the defining relation.
Even for the $W_N$ algebra, the explicit form of the defining relation
in terms of $\bar{W}^k_n$ has not been known for general $N$.
It may be easier for $q$-$W_N$.
We remark that when we study representation theory by free field
realization and BRST cohomology technique, the explicit form of the
defining relation is not necessary and what we need are $q$-Miura
transformation, screening currents and vertex operators just like
as $W_N$ case.

\noindent(\romannumeral2)~ 
Various limits.
(a) $q\rightarrow 1$: 
$\hbar^{\prime\,N}$ term of \eq{q-Miura} gives \eq{Miura}, but explicit
relation between $W^i(z)$ and $\bar{W}^k(z)$ are unknown for general $N$.
What is the meaning of $\hbar^{\prime\,\ell}$ ($\ell>N$) terms?
Are they related to conserved quantities in CFT?
(b) $q\rightarrow 0$: 
For quantum (affine) Lie algebras, $q\rightarrow 0$ limit yielded an
important notion, crystal base\cite{rKL}.
For $q$-Vir, see our another talk in this meeting.

\noindent(\romannumeral3)~ 
Kac determinant.
For $N=2$ case ($q$-Vir), we calculated and conjectured the Kac
determinant of $q$-Vir\cite{rSKAO}.
It has same structure as that of the Virasoro algebra and extra zeros
when $q$ (or $t$) is a root of unity. We guess that the Kac determinant
of $q$-$W_N$ is also so.

\noindent(\romannumeral4)~ 
Tensor product representation.
Tensor product representation of the Virasoro algebra is trivial
because the Virasoro algebra is a Lie algebra.
But $q$-Vir ($q$-$W_N$, $W_N$) satisfies quadratic relation.
So its tensor product representation is nontrivial. 

\noindent(\romannumeral5)~ 
Relation to a quantum affine Lie algebra.
Since the Virasoro algebra is obtained from affine Lie algebra 
$\hat{\cal G}$ by the Sugawara construction, we naively expect that 
$q$-Vir is obtained from quantum affine Lie algebra $U_q(\hat{\cal G})$.
In the critical level, the $q$-Virasoro Poisson algebra 
was constructed through $q$-Sugawara form\cite{rFR}.
An interesting relation between screening currents of 
$q$-$W_N$ and $U_q(\hat{\cal G})$ was pointed out in \cite{rFF}.

\noindent(\romannumeral6)~ 
Various type of $W$ algebras.
$q$-$W_N$ is $A_{N-1}$-type.
The Hamiltonian reduction technique gives us various type of $W$ algebra
from affine Lie algebras\cite{rBS}.
Can we extend this technique to the $q$-world?
$W_{1+\infty}$ algebra is obtained from $W_N$ by taking a suitable 
$N\rightarrow\infty$ limit.
What algebra is obtained from $q$-$W_N$ by $N\rightarrow\infty$ limit?
We remark that for special value of $\beta$, $q$-Vir is related to 
$W_{1+\infty}$.
Supersymmetric extension is also important
(see our another talk in this meeting).

\noindent(\romannumeral7)~ 
Macdonald operators.
The Macdonald operator can be bosonized and expressed in terms of 
$q$-Vir generators\cite{rSKAO}.
Bosonized form of the higher Macdonald operators may be expressed
by $q$-$W_N$ generators.

\noindent(\romannumeral8)~ 
Vertex operators.
Vertex operators(primary fields) are indispensable for applications
to physics. There are several proposals for vertex 
operators\cite{rAKOS,rLP2,rKa,rAKMOS,rAJMP} but at present general
definition has not been obtained.
The condition (\romannumeral1) of $q$-deformation (see footnote 1)
is of course satisfied but the problem exists in the
condition (\romannumeral2); what property should we impose on vertex
operators?

\noindent(\romannumeral9)~ 
BRST cohomology.
By using above vertex operators and screening currents, the BRST
property of the free field representation of $q$-$W_N$ can be studied.
To construct BRST charge and screened vertex operators, we should
include zero mode factor\cite{rLP2,rKa,rAKMOS,rAJMP,rJLMP}.

\noindent(\romannumeral10)~ 
Correlation functions.
By using vertex operators and screening currents, correlation functions
can be calculated\cite{rLP2,rAKMOS,rAJMP}.
In CFT correlation functions for the minimal model are characterized by
the differential equations. It is expected that correlation functions for
$q$-$W_N$ satisfy some $q$-difference equations\cite{rKa,rAKMOS}.

\noindent(\romannumeral11)~ 
Application to solvable models.
In the work\cite{rLP2}, multi-point local hight probabilities for the
ABF model in the regime \uppercase\expandafter{\romannumeral3} were
calculated, where the vertex operator of the ABF model is identified
with that of $q$-Vir. Its higher rank generalization was partially
studied in \cite{rAJMP}. 
They constructed (type \uppercase\expandafter{\romannumeral1}) vertex
operators of the RSOS model. To calculate its correlation functions,
however, the knowledge of BRST cohomology will be needed.

\noindent(\romannumeral12)~ 
Applications to other physics.
In the work\cite{rL}, the $q$-Vir generator $T(z)$ is identified with
the Zamolodchikov-Faddeev operator for the basic scalar particle in the
XYZ model. In the scaling limit XYZ model is described by the sine-Gordon
model, and the two-particle S-matrix of the basic scalar particle in the
sine-Gordon model can be obtained from the defining relation of $q$-Vir.
We remark that in CFT the meaning of $z$ in $L(z)$ is the local coordinate
of the Riemann surface, but in this case $z$ of $T(z)$ corresponds to
the rapidity of particle.

Our original motivation is to find and study the symmetry of massive
integrable models. We hope that $q$-$W_N$ will be found out to be
a useful symmetry.

\vskip5mm
\noindent{\bf Acknowledgments}\\
We would like to thank members of Nankai Institute for their
kind hospitality.
S.O. would like to thank the organizers of the meeting for giving
him the opportunity to present our results.
This work is supported in part by Grant-in-Aid for Scientific
Research from Ministry of Science and Culture. 

\section* {Appendix: Explicit examples}
\setcounter{section}{1}
\renewcommand{\thesection}{\Alph{section}}
\renewcommand{\theequation}{\Alph{section}.\arabic{equation}}
\setcounter{equation}{0}

We present explicit defining relations of $q$-$W_N$ and $q\rightarrow 1$ 
limit for $N=2,3$ cases.
To check consistency of $q\rightarrow 1$ limit (see below \eq{q-Miura-LHS}),
we need higher order terms of \eq{hin}. We set
\be
  h^i_n=\hbar'\sum_{j=1}^Nd^{ij}_n\bar{h}^j_n \quad (n\neq 0), \qquad
  d^{ij}_n=\tilde{d}^{ij}_n
  \sqrt{\frac{q^{\frac{n}{2}}-q^{-\frac{n}{2}}}{n\hbar}
        \frac{t^{\frac{n}{2}}-t^{-\frac{n}{2}}}{n\hbar\beta}}.
\ee
We remark that $d^{ij}_n$ is not uniquely determined because of 
$\sum_{i=1}^N\bar{h}^i_n=0$, and $\tilde{d}^{ij}_n$ can be chosen
as a function of $p$ only and 
$\tilde{d}^{N+1-i,N+1-j}_n=\tilde{d}^{ij}_n\Bigl|_{p\rightarrow p^{-1}}$.

\vskip 7mm
\noindent
{\large\bf $N=2$ case} ($q$-Vir)\footnote{
Since $h^1_n+p^nh^2_n=0$, it is convenient to set $h_n=p^{-\frac{n}{2}}h^1_n$.
}.  
$q$-Miura transformation is 
\be
  :\!\Bigl(p^{D_z}-\Lambda_1(z)\Bigr)
   \Bigl(p^{D_z}-\Lambda_2(p^{-1}z)\Bigr)\!:\;
  =p^{2D_z}-W^1(z)p^{D_z}+1,\quad
  W^1(z)=\Lambda_1(z)+\Lambda_2(z),
\ee
and the relation is 
\be
  f(\sfrac{w}{z})T(z)T(w)-T(w)T(z)f(\sfrac{z}{w}) \n
  =
  -\frac{(q^{\frac12}-q^{-\frac12})(t^{\frac12}-t^{-\frac12})}
        {p^{\frac12}-p^{-\frac12}}
  \Bigl(\delta(p\sfrac{w}{z})-\delta(p^{-1}\sfrac{w}{z})\Bigr),
\ee
where $T(z)=W^1(z)$ and $f(x)=f^{11}(x)$.

By multiplying $z^{\frac12}$ from the left and $z^{-\frac12}$ from the right,
$q$-Miura transformation \eq{q-Miura} becomes (see \eq{q-Miura-LHS})
\be
  (-\hbar')^2z^2
  :\!\Bigl(\alpha_0\partial_z+\partial\bar{h}^1(z)\Bigr) 
   \Bigl(\alpha_0\partial_z+\partial\bar{h}^2(z)\Bigr)\!:
  +O(\hbar^{\prime\,3})
  =
  p^{2(D_z-\frac12)}-T(z)p^{D_z-\frac12}+1.
  \label{N=2}
\ee
To check the consistency of this equation, we need $\tilde{d}^{ij}_n$, 
which is chosen as
\be
  \tilde{d}^{11}_n=\frac12p^{\frac{n}{2}}
  \sqrt{2\frac{p^{\frac{n}{2}}-p^{-\frac{n}{2}}}{p^n-p^{-n}}},\quad
  \tilde{d}^{12}_n=-\tilde{d}^{11}_n,\quad
  \tilde{d}^{3-i,3-j}_n=\tilde{d}^{ij}_n\Bigl|_{p\rightarrow p^{-1}}.
\ee
Then $T(z)$ has the following $\hbar'$ expansion,
\be
  T(z)=2+\hbar^{\prime\,2}\Bigl(z^2L(z)+\sfrac14\alpha_0^2\Bigr)
  +O(\hbar^{\prime\,4}),
\ee
where $L(z)$ is the Virasoro generator \eq{Vir}.
We remark that $T(z)$ is an even function of $\hbar'$.
Using these, $O(\hbar^{\prime\,\ell})$ terms ($\ell<2$) of 
R.H.S. of \eq{N=2} actually vanish 
and $O(\hbar^{\prime\,2})$ term gives Miura transformation of $L(z)$.

\vskip 7mm
\noindent
{\large\bf $N=3$ case} ($q$-$W_3$). 
$q$-Miura transformation is
$$
  :\!\Bigl(p^{D_z}-\Lambda_1(z)\Bigr)
   \Bigl(p^{D_z}-\Lambda_2(p^{-1}z)\Bigr)
   \Bigl(p^{D_z}-\Lambda_3(p^{-2}z)\Bigr)\!:\;
  =p^{3D_z}-W^1(z)p^{2D_z}+W^2(p^{-\frac12}z)p^{D_z}-1,
$$
\vskip -8mm
\ba
  W^1(z)
  &\!\!=\!\!&
  \Lambda_1(z)+\Lambda_2(z)+\Lambda_3(z), \\
  W^2(z)
  &\!\!=\!\!&
  :\!\Lambda_1(p^{\frac12}z)\Lambda_2(p^{-\frac12}z)\!:
  +:\!\Lambda_1(p^{\frac12}z)\Lambda_3(p^{-\frac12}z)\!:
  +:\!\Lambda_2(p^{\frac12}z)\Lambda_3(p^{-\frac12}z)\!:, \nonumber
\ea
and the relations are
\ba
  &&
  f^{11}(\sfrac{w}{z})W^1(z)W^1(w)-W^1(w)W^1(z)f^{11}(\sfrac{z}{w}) \n
  &=\!\!&
  -\frac{(q^{\frac12}-q^{-\frac12})(t^{\frac12}-t^{-\frac12})}
        {p^{\frac12}-p^{-\frac12}}
  \Bigl(\delta(p\sfrac{w}{z})W^2(p^{\frac12}w)
       -\delta(p^{-1}\sfrac{w}{z})W^2(p^{-\frac12}w)\Bigr),\n
  &&
  f^{22}(\sfrac{w}{z})W^2(z)W^2(w)-W^2(w)W^2(z)f^{22}(\sfrac{z}{w})
  \qquad\qquad \Bigl(f^{22}(x)=f^{11}(x)\Bigr) \n
  &=\!\!&
  -\frac{(q^{\frac12}-q^{-\frac12})(t^{\frac12}-t^{-\frac12})}
        {p^{\frac12}-p^{-\frac12}}
  \Bigl(\delta(p\sfrac{w}{z})W^1(p^{\frac12}w)
       -\delta(p^{-1}\sfrac{w}{z})W^1(p^{-\frac12}w)\Bigr),\\
  &&
  f^{12}(\sfrac{w}{z})W^1(z)W^2(w)-W^2(w)W^1(z)f^{21}(\sfrac{z}{w})
  \qquad\qquad \Bigl(f^{21}(x)=f^{12}(x)\Bigr) \n
  &=\!\!&
  -\frac{(q^{\frac12}-q^{-\frac12})(t^{\frac12}-t^{-\frac12})}
        {p^{\frac12}-p^{-\frac12}}
  \Bigl(\delta(p^{\frac32}\sfrac{w}{z})
       -\delta(p^{-\frac32}\sfrac{w}{z})\Bigr).\nonumber
\ea
We remark that there is no distinction between $W^1(z)$ and $W^2(z)$ in
algebraically.

By multiplying $z$ from the left and $z^{-1}$ from the right,
$q$-Miura transformation \eq{q-Miura} becomes (see \eq{q-Miura-LHS})
\ba
  &&
  (-\hbar')^3z^3
  :\!\Bigl(\alpha_0\partial_z+\partial\bar{h}^1(z)\Bigr) 
   \Bigl(\alpha_0\partial_z+\partial\bar{h}^2(z)\Bigr)
   \Bigl(\alpha_0\partial_z+\partial\bar{h}^3(z)\Bigr)\!:
  +O(\hbar^{\prime\,4}) \n
  &=\!\!&
  p^{3(D_z-1)}-W^1(z)p^{2(D_z-1)}+W^2(p^{-\frac12}z)p^{D_z-1}-1.
  \label{N=3}
\ea
We choose $\tilde{d}^{ij}_n$ as 
\ba
  \!\!\!\!&&
  \tilde{d}^{11}_n=-\tilde{d}^{13}_n=
  \frac12p^n\sqrt{2\frac{p^{\frac{n}{2}}-p^{-\frac{n}{2}}}{p^n-p^{-n}}},
  \quad
  \tilde{d}^{22}_n=-\tilde{d}^{21}_n=
  \frac12\sqrt{\frac{3}{2}
  \frac{p^n-p^{-n}}{p^{\frac{3}{2}n}-p^{-\frac{3}{2}n}}}, \n
  \!\!\!\!&&
  \tilde{d}^{12}_n=
  -\frac12p^{\frac32n}
  \sqrt{2\frac{p^{\frac{n}{2}}-p^{-\frac{n}{2}}}{p^n-p^{-n}}}
  \sqrt{3\frac{p^{\frac{n}{2}}-p^{-\frac{n}{2}}}
              {p^{\frac{3}{2}n}-p^{-\frac{3}{2}n}}},
  \quad\quad\;
  \tilde{d}^{4-i,4-j}_n=\tilde{d}^{ij}_n\Bigl|_{p\rightarrow p^{-1}}.
\ea
Then $W^1(z)$ and $W^2(z)$ have the following $\hbar'$ expansion,
\ba
  W^1(z)
  &\!\!=\!\!&
  3+\hbar^{\prime\,2}\Bigl(z^2L(z)+\alpha_0^2\Bigr) \n
  &&\quad
  +\hbar^{\prime\,3}\biggl(
  \sfrac{1}{2}z^3\Bigl(W(z)+\sfrac{1}{2}\alpha_0\partial L(z)\Bigr)
  +\sfrac{1}{4}\alpha_0z^2\Bigl(2X(z)+DX(z)\Bigr)\biggr)
  +O(\hbar^{\prime\,4}),\n
  W^2(z)
  &\!\!=\!\!&
  3+\hbar^{\prime\,2}\Bigl(z^2L(z)+\alpha_0^2\Bigr) \\
  &&\quad
  +\hbar^{\prime\,3}\biggl(
  -\sfrac{1}{2}z^3\Bigl(W(z)+\sfrac{1}{2}\alpha_0\partial L(z)\Bigr)
  +\sfrac{1}{4}\alpha_0z^2\Bigl(2X(z)+DX(z)\Bigr)\biggr)
  +O(\hbar^{\prime\,4}),\nonumber
\ea
where $L(z)=-\bar{W}^2(z)$, $W(z)=\bar{W}^3(z)$, $D=z\frac{d}{dz}$ 
and $X(z)$ is
\be
  X(z)=\sfrac12:\!\Bigl(\partial\phi^1(z)\Bigr)^2\!:
  +\alpha_0\partial^2\phi^1(z)
  -\sfrac12:\!\Bigl(\partial\phi^2(z)\Bigr)^2\!:-\alpha_0\partial^2\phi^2(z).
\ee
We remark that the combination $W(z)+\frac{1}{2}\alpha_0\partial L(z)$ 
is an primary field of the Virasoro algebra.
Using these, $O(\hbar^{\prime\,\ell})$ terms ($\ell<3$) of 
R.H.S. of \eq{N=3} actually vanish 
and $O(\hbar^{\prime\,3})$ term gives the Miura transformation of the $W_3$ 
algebra.


\end{document}